\documentclass{article}
\usepackage{latexsym}
\usepackage{epsfig}
\usepackage{amsmath,amsthm,amssymb}
\usepackage{psfrag}
\newcommand{\bfrac}[2]{\left(\frac{#1}{#2}\right)}

\def\cG{{\cal G}}
\newcommand{\brac}[1]{\left(#1\right)}
\def\cD{{\cal D}}

\def\cE{{\cal E}}

\def\E{{\bf E}}

\def\cQ{{\cal Q}}

\def\a{\alpha}

\def\d{\delta}

\def\G{\Gamma}

\def\th{\theta}

\def\l{\lambda}

\def\r{\rho}

\def\s{\sigma}
\def\S{\Sigma}

\def\OM{\Omega}

\def\Pr{\mbox{{\bf Pr}}}

\def\whp{{\bf whp}}

\def\G{\Gamma}

\newtheorem{lemma}{Lemma}
\newtheorem{theorem}{Theorem}
\newtheorem{corollary}[theorem]{Corollary}

\newcommand{\proofstart}{{\bf Proof\hspace{2em}}}

\newcommand{\proofend}{\hfill\mbox{$\Box$}}

\newcommand{\rdup}[1]{{\mbox{$ \lceil #1 \rceil $}}}
\newcommand{\rdown}[1]{{\mbox{$ \lfloor #1 \rfloor $}}}

\newtheorem{remark}{Remark}

\def\a{\alpha}

\def\Pr{\mbox{{\bf Pr}}}

\def\bE{\overline{E}}

\begin{document}

\large

\makeatletter
\title{How many random edges make a dense graph Hamiltonian?}
\author{Tom Bohman\thanks{Supported in part by NSF grant DMS-0100400}
\ \ Alan Frieze\thanks{Supported in part by NSF grant CCR-9818411.}
\ \ Ryan Martin\thanks{Supported in part by NSF VIGRE grant DMS-9819950.}  
\\ \\Department of Mathematical Sciences,\\ 
Carnegie Mellon University\\Pittsburgh PA 15213.}\date{}
\maketitle
\makeatother
\begin{abstract}
This paper investigates the number of random edges required to add to
an arbitrary dense graph in order to make the resulting graph
hamiltonian with high probability.  Adding $\Theta(n)$ random edges is
both necessary and sufficient to ensure this for all such dense
graphs.  If, however, the original graph contains no large independent
set, then many fewer random edges are required. We prove a similar
result for directed graphs.
\end{abstract}

\section{Introduction}
In the classical model of a random graph (Erd\H{o}s and R\'enyi
\cite{ER}) we add random edges to an empty graph, all at once or one
at a time, and then ask for the probability that certain structures
occur. At the present time, this model and its variants, have generated a
vast number of research papers and at least two excellent books,
Bollob\'as \cite{B0} and Janson, {\L}uczak and Ruci\'nski \cite{JLR}. 
It is also of interest to consider random graphs generated in other
ways. For example there is a well established theory of considering
random subgraphs of special graphs, such as the $n$-cube. In this
paper we take a slightly different line. We start with a graph $H$ chosen
{\em arbitrarily} from some class of graphs and then consider adding a
random set of edges $R$. We then ask if the random graph $G=H+R$ has a
certain property. This for example would model graphs which were
basically deterministically produced, but for which there is some
uncertainty about the complete structure. In any case, we feel that
there is the opportunity here for asking interesting and natural
questions.

As an example we consider the following scenario: Let $0<d<1$ be a
fixed positive constant. We let 
$\cG(n,d)$ denote the set of graphs with vertex set $[n]$ which have
minimum degree $\d\geq dn$. We choose $H$ arbitrarily from $\cG(n,d)$ 
and add a
random set of $m$ edges $R$ to create the random graph $G$. We prove
two theorems about the number of edges needed to have $G$ Hamiltonian
\whp.\footnote{A sequence of events ${\cal E}_n$ is said to 
occur ``with high probability'' (\whp) if $\lim_{n\to\infty}\Pr({\cal E}_n)=1$} Since $d\geq 1/2$ implies that $H$ itself is Hamiltonian
(Dirac's Theorem), this could be considered to be a probabilistic
generalisation of this theorem to the case where $d<1/2$. 
We henceforth assume $d<1/2$. Also, let
$$\th=\ln d^{-1}\geq .69.$$
\begin{theorem}\label{th1}\ 
Suppose $0<d\leq 1/2$ is constant, 
$H\in \cG(n,d)$. Let $G=H+R$ where
$|R|=m$ is chosen randomly from $\bE=[n]^{(2)}\setminus E(H)$. 
\begin{description}
\item[(a)] If $m\geq (30\th+13)n$ then $G$ is Hamiltonian \whp.
\item[(b)] For $d\leq 1/10$ there exist graphs $H\in \cG(n,d)$ 
such that if $m<\th n/3$
then \whp\;$G$ is not Hamiltonian.
\end{description}
\end{theorem}
So it seems that we have to add $\Theta(n)$ random edges 
in order to make $G$ Hamiltonian \whp. Since a
random member of $\cG(n,d)$ is already likely to be Hamiltonian, this is
a little disappointing. Why should we need so many edges in the worst-case?
It turns out that this is due to the existence of a large independent
set. 
Let $\a=\a(H)$ be the independence number of $H$.
\begin{theorem}\label{th2}
Suppose $H\in \cG(n,d)$ and $1\leq \a<d^2n/2$ and 
so $d>n^{-1/2}$ ($d$ need not be constant in this theorem). Let $G=H+R$ where
$|R|=m$ is chosen randomly from $\bE$. 
If
\begin{equation}\label{lb}
\frac{md^3}{\ln d^{-1}}\to\infty
\end{equation} 
then $G$ is Hamiltonian \whp.
\end{theorem}
Note that if $d$ is constant then Theorem \ref{th2} implies 
that $m\to\infty$ is sufficient. 

We have also considered a similar problem in relation to adding random arcs
to a dense digraph. For a digraph $D$ we denote its arc-set by $A(D)$. We denote
its minimum out-degree by $\d^+$
and its minimum in-degree by $\d^-$ and then we let
$\d=\min\{\d^+,\d^-\}$. Let $0<d<1$ be a fixed positive constant. 
We let 
$\cD(n,d)$ denote the set of digraphs with vertex set $[n]$ which have
$\d\geq dn$.
\begin{theorem}\label{th3}
Suppose $0<d<1/2$ is constant and 
$H\in \cD(n,d)$. Let $D=H+R$ where
$|R|=m$ is chosen randomly from $\overline{A}=[n]^2\setminus E(H)$. 
\begin{description}
\item[(a)] If $m\geq (d^{-1}(15+6\th)+
5d^{-2})n$ then $D$ is Hamiltonian \whp.
\item[(b)] For $d\leq 1/10$ there exist digraphs $H\in \cD(d)$ such 
that if $m<\th n/3$ then \whp\ $D$ is not Hamiltonian. 
\end{description}
\end{theorem}
If $\d\geq n/2$ then $H$ itself is Hamiltonian, Ghouila-Houri
\cite{GH}.

Theorem \ref{th1} is proven in the next section, Theorem \ref{th2}
is proven in Section \ref{smallalpha} and Theorem \ref{th3} is proved
in Section \ref{directed}.

\section{The worst-case for graphs}\label{wc}
We will assume from now on that $m$ is exactly $\rdup{30\th n}+13n$. We
let $R=R_1\cup R_2$ where $|R_1|=\rdup{30\th n}$. Then let $G_1=H+R_1$.

We first show that 
\begin{lemma}\label{lem1}
$G_1$ is connected \whp.
\end{lemma}
\proofstart
Let $N=\binom{n}{2}$. If $u,v\in [n]$ then either they are at distance
one or two in $H$ or
$$\Pr(dist_{G_1}(u,v)>3)\leq \left(1-\frac{|R_1|}{N}\right)^{d^2n^2}\leq 
e^{-60\th d^2n}.$$
Hence,
$$\Pr(diam(G)>3)\leq n^2e^{-60\th d^2n}=o(1).$$
\proofend

Given a longest path $P$ in a graph $\G$
with end-vertices $x_0,y$
and an edge $\{y,v\}$ where $v$ is an internal vertex of $P$, we obtain a new 
longest path $P'=x_0..vy..w$ where $w$ is the neighbor of $v$ on $P$ between
$v$ and $y$. We say that $P'$ is obtained from $P$ by a {\em rotation} with
$x_0$ {\em fixed}. 

Let $END_\G(x_0,P)$ be the set of end-vertices of longest paths of
$\G$ which can be obtained from $P$ by a sequence of rotations keeping
$x_0$ as a fixed end-vertex. 
%For each $y\in END_\G(x_0,P)$ let
%$END_\G(y,P)$ be the set of end-vertices of longest paths of
%$\G$ which can be obtained from $P$ by a sequence of rotations keeping
%$y$ as a fixed end-vertex. 
Let $END_\G(P)=\{x_0\}\cup END_\G(x_0,P)$.
Note that if $\G$ is connected and non-Hamiltonian then there is no
edge $\{x_0,y\}$ where $y\in END_\G(x_0,P)$.

It follows from P\'osa \cite{Posa} that 
\begin{equation}
\label{eq1}
|N_\G(END_\G(P))| < 2|END_\G(P)|,
\end{equation}
where for a graph $\G$ and a set $S\subseteq V(\G)$
$$N_\G(S)=\{ w\not\in S:\exists v\in S \mbox{ such that } vw\in E(\G)\}.$$
\begin{lemma}\label{lem2}
{\bf Whp} 
\begin{equation}\label{eq2}
|N_{G_1}(S)|\geq 3|S|
\end{equation}
for all $S\subseteq [n],\,|S|\leq n/5$.
\end{lemma}
\proofstart
Now $|N_H(S)|\geq 3|S|$ for all 
$S\subseteq [n],\,|S|\leq dn/3$. So,
\begin{eqnarray*}
\Pr(\exists |S|\leq n/5:\;|N_{G_1}(S)|<
3|S|)&\leq&\sum_{k=dn/3}^{n/5}\binom{n}{k}
\binom{n}{3k}\left(1-\frac{m}{N}\right)^{k(n-4k)}\\
&\leq&\sum_{k=dn/3}^{n/5}\left(\frac{n^4e^4}{27k^4}e^{-12\th}\right)^k\\
&=&o(1).
\end{eqnarray*}
\proofend

It follows from Lemma~\ref{lem2} that for any longest path \(P\) in 
a graph \( \G\) that contains \(G_1\) as a subgraph we have 
\(n/5 \le |END_{\G}(P)| \le |P|\).

Now let $R_2$ be obtained from a random sequence $e_1,e_2,\ldots$ of
edges chosen from $\bE$ {\em with replacement}.

Let $P_0$ be a longest path in $G_1$ of length $\l_0\geq dn$.
Now consider the following process: At a general stage we will 
have a path $P_i$ of length at least $\l_0+i$. We will have used a set 
$S_i\subseteq R_2$ of size $Y_i$ to go from $P_{i-1}$ to $P_i$, for 
$i\geq 1$. Here $S_1=\{e_1,e_2,\ldots,e_{Y_1}\}$, $S_2=\{e_{Y_1},
e_{Y_1+1},\ldots,e_{Y_1+Y_2}\}$ and so on. Let $Z_i=Y_1+Y_2+\cdots+Y_i$
and let $\G_i=G_1+\{e_1,e_2,\ldots,e_{Z_i}\}$. 

In order to see how the $Y_i$ are determined, let $P_i$ be a longest
path in $\G_i$ and let $END_{\G_i}(P_i)$ be as defined above and note
that by Lemma \ref{lem2}, we can assume that $|END_{\G_i}(P_i)|\geq
n/5$. We now add edges $e_{Z_i+1},e_{Z_i+2},\ldots$ in turn until we
find an edge $e_{Z_i+k}=\{a,b\}$ where $a\in END_{\G_i}(P_i)$ and
$b\in END_{\G_i}(a,P_i)$. Since $\G_i$ is connected the addition of
$\{a,b\}$ will increase the length of the longest path or close a
Hamilton cycle. We let $Y_{i+1}=k$ in this case. Finally, once we have
formed a Hamilton cycle, at stage $r$ say, we let
$Y_{r+1}=\cdots=Y_n=0$.  

Now the random variables $Y_1,Y_2,\ldots,Y_n$ are independent
random variables which are either geometric with probability of 
success at least $\frac{2}{25}$
or are zero valued. Thus 
$$\E(Z_n)\leq \frac{25n}{2}.$$
Since the variance of $Z_n$ is $O(n)$ 
it is easy to show by an application of 
Chebychev's inequality that
$Z_n\leq 13n$ \whp\ and this completes the proof of (a).
\begin{remark}
The calculations above go through quite happily for $\d(H)\geq
n^{3/4}$, say. For this degree bound the number of additional edges
required in the worst-case is $\OM(n\ln n)$. 
Since only
$\frac{1}{2}n\ln n$ random edges are required for Hamiltonicity
when we start
with the empty
graph, there is no point in considering smaller values of $d$,
unless we can improve the constant factor.
\end{remark}
(b) Let $m=cn$ for some constant $c$ 
and let $H$ be the complete bipartite graph $K_{A,B}$ where $|A|=dn$
and $|B|=(1-d)n$. Let $I$ be the set of vertices of $B$ which are not
incident with an edge in $R$. If $|I|>|A|$ then $G$ is not
Hamiltonian. Instead of choosing $m$ random edges for $R$, we choose each
possible edge independently with probability $p=\frac{2m}{(d^2+(1-d)^2)n^2}$.
(We can use monotonicity, see for example Bollob\'as \cite{B0} II.1 to 
justify this
simplification). Then
$$\E(|I|)=(1-d)n(1-p)^{(1-d)n-1}\sim
(1-d)\exp\left\{-\frac{2(1-d)m}{(d^2+(1-d)^2)n}\right\}n.$$
It follows from the Chebychev inequality that $|I|$ 
is concentrated around its mean and so $G$ will be non-Hamiltonian \whp\ 
if $c$
satisfies 
$$c<\frac{1}{2(1-d)}(d^2+(1-d)^2)\ln(d^{-1}-1).$$
This verifies (b).
\proofend
\section{Graphs with small independence number}\label{smallalpha}
{\bf Proof of Theorem \ref{th2}}\\ 
We will first show that we can decompose $H$ into a few large cycles. 
\begin{lemma}\label{decomp}
Suppose that $G$ has minimum degree $dn$ where $d\leq 1/2$ and 
that $\a(G)< d^2n/2$. Let $k_0=\rdown{\frac{2}{d}}$. Then
the vertices of $G$ can be partitioned into $\leq k_0$ 
vertex disjoint cycles.
\end{lemma}
\proofstart
Let $C_1$ be the largest cycle in $H$.  $|C_1|\geq dn+1$ 
and we now show that the graph $H\setminus C_1$ has minimum degree 
$\geq dn-\a$.

To see this, let $C_1=v_1,\ldots,v_c,v_{c+1}=v_1$.  Let $w\in 
V(H\setminus C_1)$.  Because $C_1$ is maximum sized, no such $w$ 
is adjacent to both $v_i$ and $v_{i+1}$.  Also, if $w\sim v_i$ and 
$w\sim v_j$ with $i<j$ and $v_{i-1}\sim v_{j-1}$, then
\[ w,v_j,\ldots,v_c,v_1,\ldots,v_{i-1},v_{j-1},\ldots,v_i,w \]
is a larger cycle.  So the predecessors of $N(w)$ in $C_1$ must 
form an independent set and $|N(w)\cap C_1|\leq\a$.  Similar 
arguments are to be found in \cite{CE}.

We can repeat the above argument to create disjoint cycles 
$C_1,\ldots,C_r$ where $|C_1|\geq |C_2|\geq \cdots \geq |C_r|$ and
$C_j$ is a maximum sized cycle in the graph 
$H_{j-1}=H\setminus\left(C_1\cup\cdots\cup C_{j-1}\right)$ for
$j=1,2,\ldots,r$.
Now $H_k$ has minimum degree at least $dn-k\a$ and at most 
$n-dn-1-(dn-\a+1)-\cdots-(dn-(k-1)\a+1)=n-k(dn+1-(k-1)\a/2)$ vertices.
Since $d^2n>2\a$, $H_{k_0}$, if it existed, would 
have minimum degree at least 2 and a negative number of vertices.
\proofend

Let \( C_1, \dots, C_r \) be the cycles given by Lemma~\ref{decomp}.

In order to simplify the analysis, we assume 
the edges of $R$ are chosen from
$\bE$ by including each $e\in \bE$ independently with probability
$\frac{m}{|\bE|}$. Because Hamiltonicity is a monotone property,
showing that $G$ is Hamiltonian \whp\ in this model implies the
theorem. We get a further simplification in the analysis if we 
choose these random edges in rounds: set $R=R_1\cup
R_2\cup \cdots \cup R_r$ where each edge set $R_i$ is independently chosen by
including $e\in \bE$ with probability $p$, where
$1-(1-p)^r=\frac{m}{|\bE|}$. Each $R_i$ will be used to either extend
a path or close a cycle and will only be used for one such attempt. In
this way each such attempt is independent of the previous.
To this end let $G_t=H\cup\bigcup_{s=1}^t R_t$ for
$t=0,1,\ldots,r$. Thus $G_0=H$ and $G_r=G$.

Let $e=\{x,y\}$ be an edge of $C_r$ and let $Q$ be the path
$C_r-e$. In each phase of our procedure, we have a current path $P$ with
endpoints $x,y$ together with a collection of vertex disjoint cycles
$A_1,A_2,\ldots,A_s$ which cover $V$. Initially $P=Q$, $s=r-1$ and
$A_i=C_i,\, i=1,2,\ldots,r-1$. 

Suppose a path $P$ and collection of edge disjoint cycles have
been constructed in $G_{t-1}$ (initially $t=1$).
Consider the set $Z=END_{G_{t-1}}(x,P)$ created from rotations with
$x$ as a fixed endpoint, as in Section \ref{wc}. We identify the following
possibilities:

{\bf Case 1}: There exists $z_1\in Z,\,z_2\notin P$ such that
$f=(z_1,z_2)$ is an edge of $H$. \\
Let $Q$ be the corresponding path
with endpoints $x,z_1$ which goes through $V(P)$.
Now suppose that $z_2\in C_i$ and let
$f'=(z_2,z_3)$ be an edge of $C_i$ incident with $z_2$. 
Now replace $P$ by the path $Q,f,Q'$ where $Q'=C_i-f$. This construction
reduces the number of cycles by one.

{\bf Case 2}: $|V(P)|\leq n/2$ and $z\in Z$ implies that
$N_{G_{t-1}}(z)\subseteq V(P)$.\\
It follows from (\ref{eq1}) that $|Z|\geq dn/3$. 
Now add the next set $R_t$ of random edges. 
Since $|V(P)|\leq n/2$,
the probability that no edge in $R_t$ joins $z_1\in Z$ to 
$z_2\in V\setminus V(P)$
is at most $(1-p)^{(dn/3)(n/2)}$. If there is no such edge, we fail,
otherwise we can use $(z_1,z_2)$ to proceed as in Case 1.
We also replace $t$ by $t+1$.

{\bf Case 3}: $|V(P)|> n/2$ and $z\in Z$ implies that
$N_{G_{t-1}}(z)\subseteq V(P)$.\\
Now we close $P$ to a cycle. For each $z\in Z$ let
$A_z=END_{G_{t-1}}(z,Q_z)$ where $Q_z$ is as defined in Case 1. Each
$A_z$ is of size at least $dn/3$. Add in the next set $R_t$ of random
edges. The probability that $R_t$ contains no edge of the form
$(z,z')$ where $z\in Z$ and $z'\in A_z$ is at most
$(1-p)^{d^2n^2/10}$. If there is no such edge, we fail. Otherwise, we
have constructed a cycle $C$ through the set $V(P)$ in the graph $G_t$. If 
$C$ is Hamiltonian we stop. Otherwise, we choose a remaining cycle $C'$,
distinct from $C$ and replace $P$ by $C'-e$ where $e$ is any edge of
$C'$. Now $|V(P)|<n/2$ and we can proceed to Case 1 or Case 2.

After at most $r$ executions of each of the above three cases, we
either fail or produce a Hamilton cycle. The probability of failure is
bounded by
\begin{eqnarray*}
k_0((1-p)^{(dn/3)(n/2)}+(1-p)^{d^2n^2/10})&\leq&2d^{-1}\brac{
\brac{1-\frac{m}{|\bE|}}^{\frac{dn^2}{6r}}+
\brac{1-\frac{m}{|\bE|}}^{\frac{d^2n^2}{10r}}}\\
&\leq&4d^{-1}e^{-md^3/10}\\
&=&o(1)
\end{eqnarray*}
provided (\ref{lb}) holds.
\proofend

An observation: We do not actually need the condition that
$\alpha(H)\leq d^2n/2$ to complete this proof.  The weaker
condition that $d^2n/2$ bounds the independence number of the
neighborhood of each vertex is enough.
\section{Directed graphs}\label{directed}
For a digraph $D=([n],A)$ we let $B_D$ be the bipartite graph
$([1,n],[n+1,2n],E)$ such that $B_D$
contains an edge $(x,y)$ for every arc $(x,y-n)\in A$. Perfect matchings
of $B_D$ correspond to {\em cycle covers} of $D$ i.e. sets of vertex
disjoint directed cycles which contain all vertices of $D$. 

We divide our arcs $R$ into two subsets: $R=R_1\cup R_2$, where
each $R_i$ is independently randomly chosen from $[n]^2\setminus
A(H)$. Here
$$|R_i|=\r_in\mbox{ where }\r_1=d^{-1}(15+6\th)\mbox{ and }\r_2=5d^{-2}.$$
\begin{lemma}\label{q1}
{\bf Whp} $H_1=H+R_1$ has a cycle cover $\S_1$.
\end{lemma}
\proofstart
We apply Hall's theorem to show that $B_{H_1}$ has a perfect matching
\whp. If $B_{H_1}$ does not have a perfect matching then there exists
a {\em witness} $K\subseteq [1,n],|K|\leq n/2$ (or $L\subseteq [n+1,2n],
|L|\leq n/2$) such that its neighbor set $N(K)$ in $B_{H_1}$ satisfies
$|N(K)|\leq |K|-1$ (resp. $|N(L)|\leq |L|-1$). Clearly any such witness
must be of size at least $dn$.

Since having a perfect matching is a monotone increasing property, 
we can assume that the arcs of $R_1$ are chosen independently with
$p_1=\frac{\r_1}{n}$. 

The probability that $B_{H_1}$ does not contain a perfect matching is 
therefore bounded by
$$2\sum_{k=dn}^{n/2}\binom{n}{k}\binom{n}{k-1}(1-p_1)^{k(n-k)}\leq
2\sum_{k=dn}^{n/2}\left(\frac{n^2e^2}{k^2}\cdot e^{-\r_1/2}\right)^k
=o(1).$$
\proofend

We also need to know that there are many arcs joining large sets.
For $S\subseteq [n]$ let $N^+(S)=\{t\notin S:\;\exists
s\in S$ such that $(s,t)$ is an arc of $H_1\}$. Define $N^-(S)$
similarly. 
\begin{lemma}\label{q2}
{\bf Whp}, for all disjoint $S,T\subseteq [n]$ with 
$|S|,|T|\geq dn/2$, $|N^-(S)\cap T|,|N^+(S)\cap T|\geq |T|/2$.
\end{lemma}
\proofstart
Let $\cE$ denote the event $\{\exists \mbox{ disjoint }
S,T\subseteq [n]:\;|S|,|T|\geq dn/2$ and $|N^+(S)\cap T|< |T|/2\}$.
Now fix $S,T$ with $|S|=s,|T|=t$. 
If $|N^+(S)\cap T|< |T|/2$ then there exists $T'\subseteq T$ of
size $|T|/2$ such that there are no arcs from $S$ to $T'$ in
$H_1$. The probability of this is at most
$$2^t(1-p_1)^{ts/2}\leq \brac{2e^{-\r_1s/(2n)}}^t.$$ 
Thus
\begin{eqnarray*}
\Pr(\cE)&\leq&\sum_{s=dn/2}^{(1-d/2)n}\sum_{t=dn/2}^{n-s}
\binom{n}{s}\binom{n}{t}\brac{2e^{-\r_1s/(2n)}}^t\\
&\leq&\sum_{s=dn/2}^{(1-d/2)n}\sum_{t=dn/2}^{n-s}\bfrac{ne}{s}^s
\bfrac{2ne^{1-\r_1s/(2n)}}{t}^t\\
&\leq&\sum_{s=dn/2}^{(1-d/2)n}\sum_{t=dn/2}^{n-s}
\bfrac{ne}{s}^se^{-\r_1st/(3n)}\\
&\leq&n^2e^{-\r_1dn^2/24}\\
&=&o(1).
\end{eqnarray*}
The proof for $N^-(S)\cap T$ is identical.
\proofend

\begin{corollary}\label{cor1}
$H_1$ is strongly connected \whp.
\end{corollary}
\proofstart
If $H_1$ is not strongly connected then there exists $S\subseteq [n],
|S|\leq n/2$ such that either $N^+(S)=\emptyset$ or
$N^-(S)=\emptyset$. But this would contradict Lemma \ref{q2} with
$T=\overline{S}$. 
\proofend

Assume that $H_1$ is strongly connected and satifies the condition of
Lemma \ref{q2}.

We now describe a procedure for converting the cycle cover $\S_1$ to a
Hamilton cycle. We start with an arbitrary cycle $C$ for which there
is an arc $(y,z)$ with $y\in C,z\in C'\neq C$. Such an arc must exist
because $H_1$ is strongly connected. Let $(y,x)$ be the arc of $C$
leaving $y$ and $(y',z)$ be the arc of $C'$ entering $z$. Now delete
arcs $(y,x),(y',z)$ from $\S_1$ and add the arc $(y,z)$. This yields a
path $P$ from $x$ to $y'$ plus a collection of disjoint cycles which
covers $[n]$. Call this a {\em near cycle cover} (NCC).

Given an NCC we first try to perform an {\em out path extension}. We
can do this if there is an arc $e$ joining the terminal endpoint of the
path $P$ to a vertex $v$ of one of the cycles, $C'$ say. By adding the
arc $e$ and deleting the arc of $C'$ entering $v$ we create an NCC
with one fewer cycle. Note that this construction is the same as that
of the previous paragraph, except that we do not invoke strong
connectivity. We also try to perform an analogous {\em in path extension}
by checking if there is an arc $(w,s)$ where $s$ is the start vertex of
$P$ and $w\notin P$. If such an
arc exists, we may extend the path $P$ by adding a path section at its
beginning.

We continue with these path extensions until the NCC $\S_2$ that
we have no longer admits one. Let $Q=(u_0,u_1,\ldots,u_k)$ be the
path of $\S_2$ and define the successor function $\s$ by
$\s(u_i)=u_{i+1}$ for $i<k$. Now $k\geq dn$ since $\d\geq dn$ and there are no out path
extensions available. All of $u_k$'s out neighbors are in $Q$. Let
$$S=\{u_{i-1}:\;i\leq k-dn/2,(u_k,u_i)\mbox{ is an arc of }H_1\},\
T=\{u_{k-dn/2},\ldots,u_k\}.$$
Clearly $|S|\geq dn/2$. Let $T'=N^+(S)\cap T$ so that $|T'|\geq dn/2$
by Lemma \ref{q2}. For $v\in T'$ choose $\l(v)\in S$ such that
$(\l(v),v),(u_k,\s\l(v))$ are both arcs of $H_1$. For each such $v\in
T'$ consider the path 
$$Q_v=Q+(u_k,\s\l(v))-(\l(v),\s\l(v))+(\l(v),v)-(\s^{-1}(v),v).$$
Note that $Q_v$ has the same vertex set as $Q$ and has endpoints
$u_0,\s^{-1}(v)$ --- see Figure 1.
\begin{center}
\psfrag{g}{$u_0$}
\psfrag{d}{$u_k$}
\psfrag{c}{$v$}
\psfrag{a}{$\l(v)$}
\psfrag{b}{$\s\l(v)$}
\psfrag{e}{$\s^{-1}(v)$}
\epsfig{figure=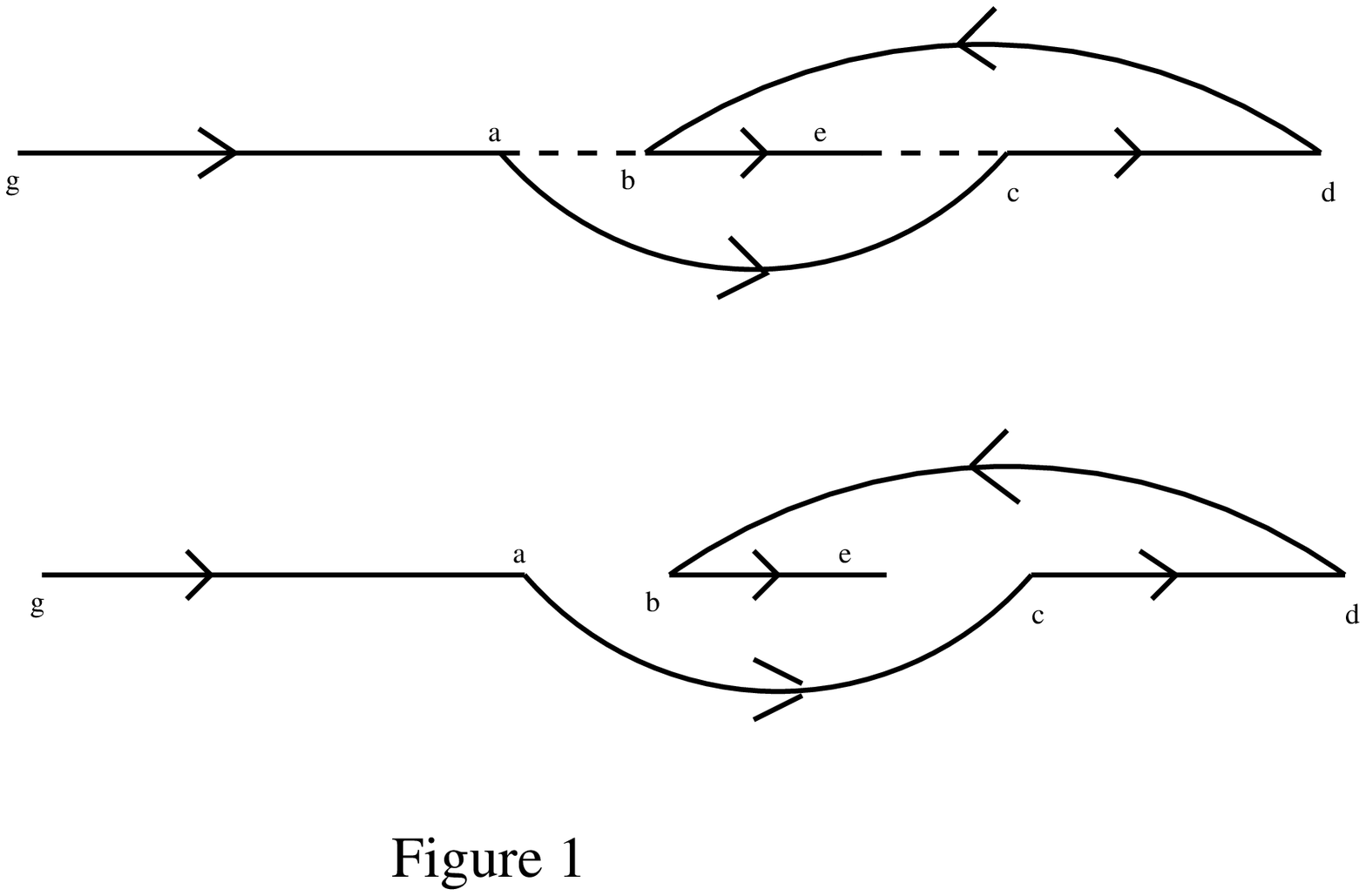,width=5in}
\end{center}
For each $v\in T'$ we see if there is an out path extension
available for $Q_v$. Suppose no such out extension exists. 
By an analogous
procedure to the creation of $Q_v$, we can, for each $v\in T'$, construct
a set $\cQ_v$ of $\geq dn/2$ paths each with a distinct start vertex
and all with the same end vertex $\s^{-1}(v)$, and all covering the
vertices of $Q$, the start vertices are distinct within $\cQ_v$ that
is. (There are no in path extensions available into
$u_0$ and we just look at the arcs that enter $u_0$).

If there is an in path extension available for a $v\in T',Q'\in \cQ_v$
then we carry it out.

Now suppose that we fail in all of these attempts at path
extension. We generate a sequence of random arcs $e_1,e_2,\ldots,$
part of $R_2$. Each $e_i$ is chosen uniformly from the arcs not in
$H_1$, with replacement. We continue until we find an arc which closes
a path in some $\cQ_v,\,v\in T'$ to a cycle $C^*$ say. Note that each
$e_i$ has probability at least $d^2/4$ of achieving this.

Now note that the sequence, starting with a cycle cover, replacing 
two cycles by a path, doing path extensions, using random arcs to
close a path to $C^*$, produces a new cycle cover with one less cycle.
Thus eventually a Hamilton cycle is produced.

The number of random edges required can be bounded by the sum
$Z=Z_1+Z_2+\cdots+Z_n$ where the $Z_i$ are independent geometric
random variables with probability of success $d^2/4$. Thus
$\E(Z)=4d^{-2}n$ and \whp\ $Z<5d^{-2}n$. We could use the Chebychev
inequality for example to prove the latter claim. Thus if we add
$5d^{-2}n$ random edges to $H_1$ then we will create a Hamilton
cycle \whp. This completes the proof of part (a) of Theorem \ref{th3}.

For part (b) we can start with $K_{A,B}$ of Theorem \ref{th1}(b) and
then replace each edge by an arc in both directions to create
$H=\vec{K}_{A,B}$. Once again we let $I$ be the set of vertices of $B$
which are not
incident with an arc in $R$. If $|I|>|A|$ then $D$ is not
Hamiltonian. Instead of choosing $m=cn$ random arcs for $R$, we choose each
possible arc independently with probability $p=\frac{m}{(d^2+(1-d)^2)n^2}$.
Then
$$\E(|I|)=(1-d)n(1-p)^{2(1-d)n-2}\sim
(1-d)\exp\left\{-\frac{2(1-d)m}{(d^2+(1-d)^2)n}\right\}n.$$
It follows from the Chebychev inequality that $|I|$ 
is concentrated around its mean and so $G$ will be non-Hamiltonian \whp\ 
if $c$
satisfies 
$$c<\frac{1}{2(1-d)}(d^2+(1-d)^2)\ln(d^{-1}-1).$$
This verifies (b).

{\bf Acknowledgement} We thank the referees for a careful reading 
which has revealed several 
small errors and led to the simplified proof of Theorem \ref{th1}(a).

\end{document}